\documentclass[11pt]{amsart}
\usepackage[utf8]{inputenc}
\usepackage{enumerate,amssymb,xcolor,comment}
\usepackage[a4paper, footskip=1cm, headheight = 16pt, top=3.7cm, bottom=3cm,  right=2.5cm,  left=2.5cm ]{geometry}
\usepackage[colorlinks,linkcolor=blue,citecolor=red]{hyperref}
\usepackage{thm-restate}

\newtheorem{question}{Question}[section]
\newtheorem{lemma}[question]{Lemma}
\newtheorem{theorem}[question]{Theorem}

\DeclareMathOperator{\cw}{cw}
\DeclareMathOperator{\rw}{rw}
\usepackage[foot]{amsaddr}

\usepackage[T1]{fontenc}

\usepackage[leqno]{amsmath}

\makeatletter
\newcommand{\leqnomode}{\tagsleft@true}
\newcommand{\reqnomode}{\tagsleft@false}
\makeatother

\usepackage{latexsym}
\usepackage{amsfonts}

\usepackage{tikz}
\usepackage{float}
\usepackage{lmodern}

\def\dd{\hbox{-}}

\usepackage{marvosym}               

\DeclareMathOperator{\tw}{tw}

\DeclareMathOperator{\cl}{cl}

\DeclareMathOperator{\comps}{cc}

\newcounter{tbox}

\newcommand{\mylongtitle}[1]{%
  \ifodd\value{page}%
    \protect\parbox{0.97\linewidth}{#1}\hfill%
  \else%
    \hfill\protect\parbox{0.97\linewidth}{#1}%
  \fi%
}

\makeatletter
\newcommand{\otherlabel}[2]{\protected@edef\@currentlabel{#2}\label{#1}}
\makeatother

\mathchardef\mh="2D

\title[Induced subgraphs and tree decompositions VI.]{Induced subgraphs and tree decompositions VI. \\ Graphs with 2-cutsets}

\author{Tara Abrishami$^{\dagger \ast}$}
\author{Maria Chudnovsky$^{\ddagger \ast}$}
\author{Sepehr Hajebi $^{\mathsection}$}
\author{Sophie Spirkl$^{\mathsection \parallel}$}
\address{$^\dagger$Department of Mathematics, University of Hamburg, Germany. This work was performed while the author was at Princeton University.}
\address{$^{\ddagger}$Princeton University, Princeton, NJ, USA}
\address{$^{\mathsection}$Department of Combinatorics and Optimization, University of Waterloo, Waterloo, Ontario, Canada}
\address{$^{\ast}$ Supported by NSF Grant DMS-1763817 and
     NSF-EPSRC Grant DMS-2120644.}
\address{$^{\parallel}$ We acknowledge the support of the Natural Sciences and Engineering Research Council of
Canada (NSERC), [funding reference number RGPIN-2020-03912]. Cette recherche a été financée
par le Conseil de recherches en sciences naturelles et en génie du Canada (CRSNG), [numéro de
référence RGPIN-2020-03912]. This project was funded in part by the Government of Ontario.}

\date {\today\\     This is the accepted manuscript; the published version appeared in Discrete Mathematics, Volume 348, Issue 1, January 2025, 114195 and is available here: \url{https://doi.org/10.1016/j.disc.2024.114195}.}

\allowdisplaybreaks

\begin{document}

\maketitle
\begin{abstract}
    This paper continues a series of papers investigating the following question: which hereditary graph classes have bounded treewidth? We call a graph {\em $t$-clean} if it does not contain as an induced subgraph the complete graph $K_t$, the complete bipartite graph $K_{t, t}$, subdivisions of a $(t \times t)$-wall, and line graphs of subdivisions of a $(t \times t)$-wall. It is known that graphs with bounded treewidth must be $t$-clean for some $t$; however, it is not true that every $t$-clean graph has bounded treewidth. In this paper, we show that three types of cutsets, namely clique cutsets, 2-cutsets, and 1-joins, interact well with treewidth and with each other, so graphs that are decomposable by these cutsets into basic classes of bounded treewidth have bounded treewidth. We apply this result to two hereditary graph classes, the class of ($ISK_4$, wheel)-free graphs and the class of graphs with no cycle with a unique chord. These classes were previously studied and decomposition theorems were obtained for both classes. Our main results are that $t$-clean ($ISK_4$, wheel)-free graphs have bounded treewidth and that $t$-clean graphs with no cycle with a unique chord have bounded treewidth. 
\end{abstract}

\section{Introduction}
All graphs in this paper are simple. A {\em tree decomposition $(T, \chi)$} of a graph $G$ consists of a tree $T$ and a map $\chi: V(T) \to 2^{V(G)}$, satisfying the following properties: 
\begin{enumerate}[(i)]
    \item For all $v \in V(G)$, there exists $t \in V(T)$ with $v \in \chi(t)$; 
    
    \item For all $uv \in E(G)$, there exists $t \in V(T)$ with $u, v \in \chi(t)$; and 
    
    \item For all $v \in V(G)$, the set $\{t \in V(T): v \in \chi(t)\}$ induces a connected subtree of $T$.
\end{enumerate}
The {\em width} of a tree decomposition $(T, \chi)$ is equal to $\max_{t \in V(T)} |\chi(t)| - 1$. The {\em treewidth} of a graph $G$ is the minimum width of a tree decomposition of $G$. Treewidth is an important graph parameter for a number of reasons. Graphs with bounded treewidth have nice structural and algorithmic properties; in particular, many \textsf{NP}-hard algorithmic problems can be solved in polynomial time in graphs with bounded treewidth. In their celebrated Grid Minor Theorem, Robertson and Seymour provided a complete characterization of graphs with bounded treewidth with respect to the subgraph relation. In particular, they proved that every graph of large treewidth contains a subdivision of a graph called a {\em wall} as a subgraph, and this subgraph has big treewidth.  See \cite{wallpaper}, page 378 for a full definition of a wall. 

\begin{lemma}[\cite{GMV}]
\label{lemma:gridthm}
There is a function $f : \mathbb{N} \to \mathbb{N}$ such that every graph of treewidth at least
$f(k)$ contains a subdivision of the $(k \times k)$-wall as a subgraph. \end{lemma}
An {\em induced subgraph} of a graph $G$ is a graph $H$ with $V(H) \subseteq V(G)$ and for all $u, v \in V(H)$, $uv \in E(H)$ if and only if $uv \in E(G)$. For a graph $H$, a graph $G$ is {\em $H$-free} if $G$ does not contain an induced subgraph isomorphic to $H$. For a set of graphs $\mathcal{H}$, a graph $G$ is {\em $\mathcal{H}$-free} if $G$ is $H$-free for all $H \in \mathcal{H}$. Graph classes defined by forbidden induced subgraphs (also known as {\em hereditary graph classes}) are important areas of interest and research in structural graph theory, but unlike in the case of subgraphs, the relationship between induced subgraphs and treewidth is not well understood. There are four well-known graph classes that are induced subgraph obstructions to bounded treewidth: complete graphs, complete bipartite graphs, subdivided walls, and line graphs of subdivided walls. Any hereditary graph class with bounded treewidth must exclude a fixed complete graph, a fixed complete bipartite graph, subdivisions of a fixed wall, and line graphs of subdivisions of a fixed wall.  Recently, Korhonen \cite{Korhonen} proved that in the case of bounded maximum degree, these are the only induced subgraph obstructions to bounded treewidth. Therefore, the focus now is to determine which graph classes of unbounded degree have bounded treewidth. Several previous papers in this series have proven that certain hereditary graph classes of unbounded degree have bounded treewidth; see \cite{onenbr, pyramiddiamond}. Graph classes in which treewidth is bounded by a logarithmic function of the number of vertices have also been studied (\cite{logpaper, stephan}).

 A graph $G$ is {\em $t$-clean} if $G$ does not contain a subdivision of the $(t \times t)$-wall, the line graph of a subdivision of the $(t \times t)$-wall, $K_t$, or $K_{t, t}$ as an induced subgraph. Note that any graph of bounded treewidth must be $t$-clean for some $t \geq 0$. Let $\mathcal{F}$ be a hereditary graph class. We say that $\mathcal{F}$ is {\em special} if every $t$-clean graph in $\mathcal{F}$ has bounded treewidth. A graph $G$ is {\em $ISK_4$-free} if $G$ does not contain a subdivision of $K_4$ as an induced subgraph. A {\em wheel $(H, v)$} consists of a hole $H$ and a vertex $v$ such that $v$ has at least three neighbors in $H$. A cycle $C$ has a {\em unique chord} if it has exactly one chord (see also the definition below). A graph $G$ does not contain a cycle with a unique chord if every cycle of $G$ has either no chords or at least two chords. Previously, ($ISK_4$, wheel)-free graphs and the class of graphs that do not contain a cycle with a unique chord were studied and structure theorems obtained for both classes (\cite{ISK4, Vuskovic-Trotignon}). In this paper, we use these structure theorems and results about treewidth and certain cutsets to prove that ($ISK_4$, wheel)-free graphs are special and that graphs that do not contain a cycle with a unique chord are special. 
 
 Previous results in this series use ``central bag methods'' to decompose graphs along well-chosen cutsets such as clique cutsets and star cutsets. Roughly, the central bag method allows us to show that the treewidth of a ``less complicated'' induced subgraph $\beta$ of $G$ is bounded, and that $\tw(G) \leq f(\tw(\beta))$ for some function $f$.  In this paper, we also decompose along well-chosen cutsets that interact well with treewidth, and so we are likewise able to simplify the problem of bounding the treewidth of $G$ to the problem of bounding the treewidth of an induced subgraph of $G$. It turns out that the different cutsets present in these graph classes (clique cutsets, 2-cutsets, and 1-joins) also interact well with each other, so the full strength of the central bag method is unnecessary. Instead, in this case, relating the treewidth of $G$ to the treewidth of the ``less complicated'' induced subgraph is easy. The results in this paper thus follow relatively straightforwardly from the properties of these cutsets and from existing structure theorems about the graph classes. Our results about the cutsets, in particular Theorem \ref{thm:prime}, could be applied to other graph classes that are decomposable along 2-cutsets. 
 
This paper is organized as follows. In Section \ref{sec:cutsets}, we describe the relationship between 2-cutsets and treewidth. In Section \ref{sec:ISK4}, we prove that ($ISK_4$, wheel)-free graphs are special. In Section \ref{sec:unique-chords}, we prove that graphs with no cycle with a unique chord are special. 
\subsection{Definitions} 
In this paper, we use induced subgraphs and their vertex sets interchangeably. A {\em path} is a graph $P$ with vertex set $\{p_1, \hdots, p_k\}$ and edge set $\{p_1p_2, \hdots, p_{k-1}p_k\}$. We write $P = p_1 \dd \cdots \dd p_k$. The {\em length} of $P$ is the number of vertices of $P$. The vertices $p_1$ and $p_k$ are the {\em ends} of $P$. The {\em interior of $P$}, denoted $P^*$, is $P \setminus \{p_1, p_k\}$. A {\em cycle} is a graph $C$ with vertex set $\{c_1, \hdots, c_k\}$ and edge set that includes $\{c_1c_2, \hdots, c_{k-1}c_k, c_kc_1\}$. We write $C = c_1 \dd \cdots \dd c_k \dd c_1$. An edge of a cycle that is not listed above is called a {\em chord}. (Note that this is a nonstandard way to define a cycle). A cycle with no chords is called a {\em hole}. Let $G$ be a graph. A {\em clique of $G$} is a subset $C \subseteq V(G)$ such that every pair of vertices in $C$ are adjacent. Two disjoint sets $X, Y \subseteq V(G)$ are {\em complete} ({\em anticomplete}) if $xy \in E(G)$ ($xy \not \in E(G)$) for all $x \in X$, $y \in Y$. The {\em neighborhood} of a vertex $v$, denoted $N(v)$, is defined as $N(v) = \{u \in V(G) : uv \in E(G)\}$. If $uv \in E(G)$, then the graph $G'$ formed by {\em subdividing edge $uv$} is given by $V(G') = V(G) \cup \{w\}$, $E(G') = (E(G) \setminus \{uv\}) \cup \{uw, wv\}$.

\section{2-cutsets and treewidth}
\label{sec:cutsets}

In this section, we describe several results about how clique cutsets and 2-cutsets interact with treewidth. Let $G$ be a graph. A {\em clique cutset} of $G$ is a clique $C$ of $G$ such that $G \setminus C$ has more connected components than $G$. A {\em 2-cutset} is a pair of vertices $u, v \in V(G)$ such that $G \setminus \{u, v\}$ is not connected. A 2-cutset $\{u, v\}$ is {\em proper} if $G \setminus \{u, v\}$ can be partitioned into two parts $(X, Y)$ such that $X$ is anticomplete to $Y$ and $|X|, |Y| \geq 2$. Let $G$ be a graph and let $\{u, v\}$ be a proper 2-cutset of $G$. Let $C$ be a connected component of $G \setminus \{u, v\}$. The {\em closure} of $C$, denoted $\cl(C)$, is obtained from the graph induced by $C \cup \{u, v\}$ by adding edge $uv$. The {\em strong closure} of $C$, denoted $\cl^*(C)$, is the graph formed from $\cl(C)$ by subdividing the edge $uv$. 

First, we remark that clique cutsets do not affect
treewidth.  For $X \subseteq V(G)$, we let $\comps(G \setminus X)$ denote the set of connected components of $G \setminus X$. 
\begin{lemma}[\cite{clique-cutsets-tw}]
Let $G$ be a graph. Then, the treewidth of $G$ is equal to the maximum treewidth
over all induced subgraphs of $G$ with no clique cutset.
\label{lemma:clique-cutsets-tw}
\end{lemma}
Next, we state a similar result for 2-cutsets. A connected graph $G$ has a {\em 1-cutset} if there exists $v \in V(G)$ such that $G \setminus \{v\}$ has at least two connected components. Note that a 1-cutset is a clique cutset. 

\begin{lemma}[\cite{clique-cutsets-tw}]
Let $G$ be a graph with no 1-cutset and let $\{u, v\}$ be a 2-cutset of $G$. Then, $\tw(G) = \max_{C \in \comps(G \setminus\{u, v\})} \tw(\cl(C))$. 
\label{lemma:2-cutset-subdivided-tw-1}
\end{lemma}

The following lemma states that subdividing an edge does not increase treewidth. 

\begin{lemma}
Suppose $G$ is a graph with $\tw(G) =k$. Let $uv$ be an edge of $G$, and let $G'$ be the graph formed by subdividing edge $uv$ exactly once. Then, $\tw(G') = k$. 
\end{lemma}
\begin{proof}
If $G$ is a forest, then $G'$ is a forest and the result holds, so we may assume $k \geq 2$. Let $w$ be the unique vertex of $V(G') \setminus V(G)$. Let $(T, \chi)$ be a tree decomposition of $G$ of width $k$. Let $t \in V(T)$ be such that $\{u, v\} \subseteq \chi(t)$. Let $T'$ be the tree formed from $T$ by adding vertex $t'$ adjacent to vertex $t$, and let $\chi'$ be the map given by $\chi'(t) = \chi(t)$ if $t \neq t'$, and $\chi'(t') = \{u, w, v\}$. Now, $(T', \chi')$ is a tree decomposition of $G'$ with width equal to $\max(2, k) = k$. This completes the proof.
\end{proof}
Now, we have the following corollary: 
\begin{lemma}
\label{lemma:2-cutset-subdivided-tw}
Let $G$ be a graph with no 1-cutset and let $\{u, v\}$ be a 2-cutset of $G$.  Then, $\tw(G) = \max_{C \in \comps({G \setminus \{u, v\})}} \tw(\cl^*(C))$. 
\end{lemma}

Finally, we state the main result of this section. Let $\mathcal{F}$ be a hereditary graph class. We say that $\mathcal{F}$ is {\em 2-cutset-safe} if for every graph $G \in \mathcal{F}$ with no 1-cutset, for every 2-cutset $\{u, v\}$ of $G$, and every component $C$ of $G \setminus \{u, v\}$, either $\cl(C) \in \mathcal{F}$ or $\cl^*(C) \in \mathcal{F}$.

A graph $G$ is {\em prime} if $G$ has no clique cutset and no proper 2-cutset. Let $\mathcal{F}$ be a 2-cutset-safe graph class and let $G \in \mathcal{F}$. We construct an {\em $\mathcal{F}$ prime decomposition} of $G$ as follows. Let $G_0 = G$. If $G$ is prime, the $\mathcal{F}$ prime decomposition of $G$ is $(G_0)$. For $i \geq 1$, the sequence is defined as follows. If $G_{i-1}$ is prime, the sequence ends. If $i$ is odd, then $G_i$ is an induced subgraph of $G_{i-1}$ with no clique cutset, and subject to that, of maximum treewidth. If $i$ is even, then $G_{i-1}$ has no clique cutset but is not prime, so $G_{i-1}$ has a proper 2-cutset $\{u, v\}$.  Let $C_1, \hdots, C_m$ be the components of $G_{i-1} \setminus \{u, v\}$, and let $G_i \in \mathcal{F} \cap \{\cl(C_1), \hdots, \cl(C_m), \cl^*(C_1), \hdots, \cl^*(C_m)\}$ with maximum treewidth.  Now, $(G_0, \hdots, G_k)$ is an $\mathcal{F}$ prime decomposition of $G$, where $G_k$ is prime. We say that $G_k$ is the {\em prime base} of the decomposition.

\begin{theorem}
Let $t > 0$ and let $\mathcal{F}$ be a 2-cutset-safe hereditary graph class. Let $G \in \mathcal{F}$ and let $H$ be the prime base of an $\mathcal{F}$ prime decomposition of $G$.  Then, $\tw(G) = \tw(H)$. Further, let $\mathcal{J} \subseteq \mathcal{F}$ be the set of all prime graphs of $\mathcal{F}$. Then, $\tw(G) \leq \max_{J \in \mathcal{J}} \tw(J)$. 
\label{thm:prime}
\end{theorem}
\begin{proof}
Let $(G_0, \hdots, G_k)$ be an $\mathcal{F}$ prime decomposition of $G$ with $G_k = H$. For $1 \leq i \leq k$, if $i$ is odd, $\tw(G_i) = \tw(G_{i-1})$ by Lemma \ref{lemma:clique-cutsets-tw}, and if $i$ is even, $\tw(G_i) = \tw(G_{i-1})$ by Lemmas \ref{lemma:2-cutset-subdivided-tw-1} and \ref{lemma:2-cutset-subdivided-tw}. It follows that $\tw(G) = \tw(H)$. 

Since $H \in \mathcal{J}$ and $\tw(G) = \tw(H)$, it follows that $\tw(G) \leq \max_{J \in \mathcal{J}} \tw(J)$. 
\end{proof}

\section{($ISK_4$, wheel)-free graphs} 
\label{sec:ISK4}
 In this section, we use Theorem \ref{thm:prime} to prove that $t$-clean ($ISK_4$, wheel)-free graphs have bounded treewidth. We first state a structure theorem for ($ISK_4$, wheel)-free graphs. Following the notation of \cite{ISK4}, we define a {\em long rich square} to be a graph $G$ containing a four-vertex hole $S = u_1 \dd u_2 \dd u_3 \dd u_4$ as an induced subgraph such that every connected component $C$ of $G \setminus S$ is a path $P = p_1 \dd \cdots \dd p_k$ with $k > 1$, $P^*$ is anticomplete to $S$, and either $N(p_1) = \{u_1, u_2\}$ and $N(p_k) = \{u_3, u_4\}$, or $N(p_1) = \{u_1, u_4\}$ and $N(p_k) = \{u_2, u_3\}$. We call $S$ the {\em central square} of the long rich square $G$. A graph $G$ is {\em chordless} if every cycle of $G$ has no chord.
\begin{theorem}[\cite{ISK4}]
\label{thm:ISK4}
Let G be an ($ISK_4$, wheel)-free graph. Then either:
\begin{itemize}
    \item $G$ is series-parallel; 
    \item $G$ is the line graph of a chordless graph with maximum degree at most three;
    \item $G$ is a complete bipartite graph; 
    \item $G$ is a long rich square;
    \item $G$ has a clique-cutset or a proper $2$-cutset. 
\end{itemize}
\end{theorem}

Note that in \cite{ISK4}, a stronger definition of proper 2-cutset is given, but it also satisfies the definition given in this paper. Next, we state a useful lemma relating the treewidth of a graph to the treewidth of its line graph. 

\begin{lemma}[\cite{harvey-wood}]
\label{lemma:tw-linegraphs}
For all graphs $G$ with maximum degree $\Delta(G)$, it holds that $$\tw(L(G)) \leq (\tw(G) + 1)\Delta(G) - 1.$$
\end{lemma}
Now, we prove that ($ISK_4$, wheel)-free graphs are 2-cutset-safe. 
\begin{lemma}
\label{lemma:ISK4-2-cutset-safe}
The class of ($ISK_4$, wheel)-free graphs is 2-cutset-safe. 
\end{lemma}
\begin{proof}
Let $G$ be an ($ISK_4$, wheel)-free graph with no 1-cutset, let $\{u, v\}$ be a 2-cutset of $G$, let $C$ be a connected component of $G \setminus \{u, v\}$. Let $C'$ be a connected component of $G \setminus \{u, v\}$ such that $C' \neq C$. Since $G$ does not have a 1-cutset, it follows that $u$ and $v$ both have neighbors in $C'$. Let $P$ be a path from $u$ to $v$ in $C'$. Since $P \setminus \{u, v\} \subseteq C'$, it follows that $P \setminus \{u, v\}$ is anticomplete to $C$. 

We claim that $\cl^*(C)$ is ($ISK_4$, wheel)-free. Let $w$ be the unique vertex of $\cl^*(C) \setminus G$. Suppose $\cl^*(C)$ contains an $ISK_4$ or a wheel, say $H$. Since $C$ is an induced subgraph of $G$ and $G$ is ($ISK_4$, wheel)-free, it follows that $w \in H$. Note that $w$ has degree 2, so $w$ is a vertex of a subdivided path. Now, the subgraph $H'$ of $G$ induced by $(H \setminus \{w\}) \cup P$ is an induced subdivision of $H$. If $H$ is an $ISK_4$, then so is $H'$; if $H$ is a wheel, then $H'$ is a wheel or an $ISK_4$, a contradiction. This completes the proof.
\end{proof}
Finally, we state the main result of this section. 
\begin{theorem}
The class of ($ISK_4$, wheel)-free graphs is special. 
\end{theorem}
\begin{proof}
Let $G$ be a $t$-clean ($ISK4$, wheel)-free graph. We show that $\tw(G) \leq 3f(t) + 2$, where $f$ is the function from Lemma \ref{lemma:gridthm}. By Lemma \ref{lemma:ISK4-2-cutset-safe}, the class of ($ISK_4$, wheel)-free graphs is 2-cutset-safe. Therefore, by Theorem \ref{thm:prime}, we may assume that $G$ does not have a clique cutset or a proper 2-cutset. Now we apply Theorem \ref{thm:ISK4}. If $G$ is series-parallel, then $\tw(G) \leq 2$ (\cite{series-parallel}). Suppose $G$ is the line graph of a chordless graph $H$ with maximum degree at most three. Then, since $G$ does not contain the line graph of a subdivision of the $(t \times t)$-wall as an induced subgraph, it follows that $H$ does not contain a subdivision of the $(t \times t)$-wall as a subgraph (note that edges become vertices in the line graph). By Lemma \ref{lemma:gridthm}, $\tw(H) \leq f(t)$, and thus by Lemma \ref{lemma:tw-linegraphs}, $\tw(G) \leq 3f(t) + 2$. If $G$ is a complete bipartite graph, then $\tw(G) \leq t$ (since $G$ is $t$-clean). Suppose $G$ is a long rich square with central square $S$. Then, $G \setminus S$ is a forest, so $\tw(G) \leq 5$.
This completes the proof. 
\end{proof}

\section{1-joins and unique-chord-free graphs} 
\label{sec:unique-chords}

In this section, we prove that $t$-clean graphs with no cycle with a unique chord have bounded treewidth. We use a structure theorem proven in \cite{Vuskovic-Trotignon}.  A {\em 1-join} of $G$ is a partition of $V(G)$ into two sets $X$ and $Y$ with $|X|, |Y| \geq 2$ such that there exist non-empty $A \subseteq X$ and $B \subseteq Y$ with $A$ complete to $B$, $X$ anticomplete to $Y \setminus B$, and $Y$ anticomplete to $X \setminus A$. We denote this 1-join by $(X, Y, A, B)$. In this section, we use two graph parameters similar to treewidth, called {\em rankwidth} and {\em cliquewidth}. We omit full definitions of these parameters here, since we are only concerned with how they relate to treewidth. The definitions can be found in \cite{clique-rank}, section 3 (cliquewidth) and section 6 (rankwidth). We denote the rankwidth and the cliquewidth of a graph $G$ by $\rw(G)$ and $\cw(G)$, respectively. 

First, we state a result from \cite{rw} about 1-joins and rankwidth. 
\begin{lemma}[\cite{rw}]
\label{lemma:1-join-rw}
Let $G$ be a graph and let $(X, Y, A, B)$ be a 1-join of $G$. Let $a \in A$ and let $b \in B$. Then, $\rw(G) = \max(\rw(X \cup \{b\}),  \rw(Y \cup \{a\}))$.  
\end{lemma}

The following results state relationships between rankwidth, cliquewidth, and treewidth. 
\begin{lemma}[\cite{rw}]
\label{lemma:rw}
Let $G$ be a graph. Then, $\rw(G) \leq \cw(G) \leq 2^{\rw(G) + 1} - 1$.
\end{lemma}
\begin{lemma}[\cite{courcelle}]
\label{lemma:cw<tw}
Let $G$ be a graph. Then, $\cw(G) \leq 2^{\tw(G)+1} + 1$. 
\end{lemma}
\begin{lemma}[\cite{gurski}]
\label{lemma:rw-2}
Let $t > 0$ and let $G$ be a graph with no $K_{t, t}$ subgraph. Then, $\tw(G) \leq 3(t-1)\cw(G) - 1$. 
\end{lemma}
Let $R(t, t)$ denote the Ramsey number, that is, every graph with $|V(G)| \geq R(t, t)$ has either a clique or a stable set of size $t$. Now, we have a corollary of the previous lemma. 
\begin{lemma}
\label{lemma:rw-corollary}
Let $t > 0$ and let $G$ be ($K_t, K_{t,t}$)-free. Then $\tw(G) \leq 3(R(t, t) - 1)2^{\rw(G) + 1} - 1$. 
\end{lemma}
\begin{proof}
We claim that $G$ has no $K_{s, s}$ subgraph, where $s = R(t, t)$. Suppose $G$ has a $K_{s, s}$ subgraph with bipartition $(X, Y)$. Since $|X| = R(t, t)$, it follows that $X$ contains a clique or an independent set of size $t$. Since $G$ is $K_t$-free, it follows that $X$ contains an independent set of size $t$, say $A$. By symmetry, $Y$ contains an independent set of size $t$, say $B$. But $(A, B)$ is the bipartition of an induced $K_{t, t}$ of $G$, a contradiction. Now, by Lemma \ref{lemma:rw-2}, $\tw(G) \leq 3(R(t, t) - 1) \cw(G) - 1$, and by Lemma \ref{lemma:rw}, $\cw(G) \leq 2^{\rw(G) + 1} - 1$, and the result follows. 
\end{proof}

Next, we state a structure theorem from \cite{Vuskovic-Trotignon} for graphs that do not contain a cycle with a unique chord. We begin with several definitions. The {\em Petersen graph} and the {\em Heawood graph} are two well-studied graphs on ten and fourteen vertices, respectively. A graph $G$ is {\em strongly 2-bipartite} if $G$ is $C_4$-free, bipartite, and there exists a bipartition $(X, Y)$ of $G$ such that every $x \in X$ has degree exactly 2 and every $y \in Y$ has degree at least 3. A graph $G$ has a {\em proper 1-join} if it has a 1-join $(X, Y, A, B)$ where $A$ and $B$ are stable sets of $G$ of size at least 2.

In \cite{Vuskovic-Trotignon}, the following theorem is proven:

\begin{theorem}[\cite{Vuskovic-Trotignon}]
\label{thm:unique-chord}
Let $G$ be a connected graph that does not contain a cycle with a unique chord. Then either:
\begin{itemize}
\item $G$ is strongly 2-bipartite;
\item $G$ is a hole
of length at least 7; 
\item $G$ is a complete graph; 
\item $G$ is an induced subgraph of the
Petersen graph or the Heawood graph; 
\item $G$ has a 1-cutset, a proper 2-cutset, or a
proper 1-join.
\end{itemize}
\end{theorem}
Note that in \cite{Vuskovic-Trotignon}, a stronger definition of proper 2-cutset is given, but it also satisfies the definition given in this paper. 
In the following lemma, we show that if $G$ is strongly 2-bipartite and $t$-clean, then the treewidth of $G$ is bounded. 
\begin{lemma}
\label{lemma:2-bipartite}
Let $G$ be a $t$-clean strongly 2-bipartite graph. Then, there exists a function $f$ such that $\tw(G) \leq f(t)$.  
\end{lemma}
\begin{proof}
For a graph $H$, let $S_H \subseteq V(H)$ denote the set of vertices of $H$ of degree at least three. By Theorem 7.1 of \cite{onenbr}, there exists a constant $w(t, \Delta)$ such that for every graph $H$, if $H$ has treewidth at least $w(t, \Delta)$ and $S_H$ has maximum degree at most $\Delta$, then $H$ contains a subdivision of the $(t \times t)$-wall or the line graph of a subdivision of the $(t \times t)$-wall as an induced subgraph. We apply this theorem to $G$. By definition of strongly 2-bipartite, it follows that $S_G$ is a partite set of a bipartition of $G$, and thus has maximum degree 0. Let $f(t) = w(t, 0)$, where $w$ is as in Theorem 7.1 of \cite{onenbr}. Now, since $G$ is $t$-clean, it follows that $\tw(G) \leq f(t)$.
\end{proof}
Next, we prove that 1-joins interact well with 1-cutsets and proper 2-cutsets. 
\begin{lemma}
\label{lemma:clique-cutset-1join}
Let $G$ be a graph with no proper 1-join, let $u$ be a 1-cutset of $G$, let $C$ be a component of $G \setminus \{u\}$, and let $\overline{C} = C \cup \{u\}$.  Then, $\overline{C}$ has no proper 1-join.
\end{lemma}
\begin{proof}
Suppose that $\overline{C}$ has a proper 1-join $(X, Y, A, B)$. By symmetry, we may assume $u \in X$. Let $C_1, \hdots, C_m$ be the components of $G \setminus \{u\}$ other than $C$. For each $C_i$, it holds that $N(C_i) = \{u\}$. Now, $(X \cup C_1 \cup \cdots \cup C_m, Y, A, B)$ is a proper 1-join of $G$, a contradiction. 
\end{proof}

\begin{lemma}
\label{lemma:2-cutset-1join}
Let $G$ be a graph with no proper 1-join, let $\{u, v\}$ be a proper 2-cutset of $G$, and let $C$ be a component of $G \setminus \{u, v\}$. Then, $\cl^*(C)$ has no proper 1-join.
\end{lemma}
\begin{proof}
Suppose $\cl^*(C)$ has a proper 1-join $(X, Y, A, B)$. Let $w$ be the vertex formed by subdividing edge $uv$. First, we claim that up to symmetry, $\{u, v\} \subseteq X$. Suppose for a contradiction that $u \in X$ and $v \in Y$. Since $u$ and $v$ are non-adjacent, $A$ is complete to $B$, and there are no other edges from $X$ to $Y$, we may assume that $u \in A$, $w \in B$, and $v \in Y \setminus B$. But $A$ is a stable set of size at least 2, so it follows that $|N(w)| \geq 3$, a contradiction. Therefore, we may assume that $\{u, v\} \subseteq X$. Let $C_1, \hdots, C_m$ be the components of $G \setminus \{u, v\}$ different from $C$. For each $C_i$ it holds that $N(C_i) \subseteq \{u, v\}$. Now, $(X \cup C_1 \cup \cdots \cup C_m, Y, A, B)$ is a proper 1-join of $G$, a contradiction. 
\end{proof}

Now, we prove that the class of graphs with no cycle with a unique chord is 2-cutset-safe. 
\begin{lemma}
\label{lemma:unique-chord-2-cutset-safe}
The class of graphs with no cycle with a unique chord is 2-cutset-safe. 
\end{lemma}
\begin{proof}
Let $G$ be a graph with no cycle with a unique chord and assume that $G$ has no 1-cutset. Let $\{u, v\}$ be a 2-cutset of $G$, let $C$ be a connected component of $G \setminus \{u, v\}$, and let $w$ be the unique vertex of $\cl^*(C) \setminus G$. Let $C'$ be a connected component of $G \setminus \{u, v\}$ such that $C' \neq C$. Since $G$ does not have a 1-cutset, it follows that $u$ and $v$ both have neighbors in $C'$. Let $P$ be a path from $u$ to $v$ in $C'$. Since $P \setminus \{u, v\} \subseteq C'$, it follows that $P \setminus \{u, v\}$ is anticomplete to $C$. 

Suppose that $\cl^*(C)$ contains a cycle with a unique chord, say $H = h_1 \dd \cdots \dd h_k \dd h_1$ with chord $h_ih_j$, where $1 \leq i < j - 1 \leq k$. Since $C$ is an induced subgraph of $G$ and $G$ does not contain a cycle with a unique chord, it follows that $w \in H$, and so $u, v \in H$. Because $w$ has degree 2, it follows that $w \not \in \{h_i, h_j\}$. Let $H'$ be the cycle of $G$ given by $(H \setminus \{w\}) \cup P$. Now, $H'$ is a cycle with a unique chord, a contradiction. It follows that the class of graphs with no cycle with a unique chord is 2-cutset-safe. 
\end{proof}
Finally we state and prove the main theorem of this section. 

\begin{theorem}
The class of graphs that do not contain a cycle with a unique chord is special. 
\end{theorem}
\begin{proof}
Let $G$ be a $t$-clean graph that does not contain a cycle with a unique chord. We show that $\tw(G) \leq 3(R(t, t) -1)2^{2^{s+1}+2} - 1$ for $s = \max(f(t), 14)$ (where $f(t)$ is as in Lemma \ref{lemma:2-bipartite}). Suppose that $G$ has a proper 1-join $(X, Y, A, B)$. Let $a \in A$ and $b \in B$. Let $G_1 = X \cup \{b\}$ if $\rw(X \cup \{b\}) \geq \rw(Y \cup \{a\})$, and let $G_1 = Y \cup \{a\}$ otherwise. For $i \geq 1$, if $G_i$ has a proper 1-join, define $G_{i+1}$ similarly with respect to $G_{i}$, until $G_k$ has no proper 1-join. Now, it follows from Lemma \ref{lemma:1-join-rw} that $\rw(G) \leq \rw(G_k)$. 

Next, we determine the treewidth of $G_k$. By Lemma \ref{lemma:unique-chord-2-cutset-safe}, $G_k$ is a member of a 2-cutset-safe graph class. Let $H_0 = G_k$ and let $(H_0, \hdots, H_m)$ be an $\mathcal{F}$ prime decomposition of $G_k$, 
where $\mathcal{F}$ is the class of graphs that do not contain a cycle with a unique chord. For $1 \leq i \leq m$, if $i$ is odd, then $H_i$ has no proper 1-join by Lemma \ref{lemma:clique-cutset-1join}, and if $i$ is even, then $H_i$ has no proper 1-join by Lemma \ref{lemma:2-cutset-1join}. Therefore, $H_m$ has no proper 1-join, and by Theorem \ref{thm:prime}, $\tw(G_k) = \tw(H_m)$. 

Let $H = H_m$. Now, we apply Theorem \ref{thm:unique-chord}. If $H$ is a clique, then $\tw(H) \leq t$. If $H$ is strongly 2-bipartite, then by Lemma \ref{lemma:2-bipartite}, $\tw(H) \leq f(t)$. If $H$ is a hole, then $\tw(H) \leq 2$. If $H$ is an induced subgraph of the Petersen graph, then $\tw(G) \leq 10$ (the number of vertices in the Petersen graph). If $H$ is an induced subgraph of the Heawood graph, then $\tw(G) \leq 14$ (the number of vertices in the Heawood graph). Therefore, $\tw(H) \leq \max(f(t), 14)$, and thus $\tw(G_k) \leq \max(f(t), 14)$. 

Let $s = \max(f(t), 14)$. By Lemma \ref{lemma:cw<tw}, $\cw(G_k) \leq 2^{s+1} + 1$, and by Lemma \ref{lemma:rw}, $\rw(G_k) \leq 2^{s+1} + 1$. By Lemma \ref{lemma:1-join-rw}, it follows that $\rw(G) \leq 2^{s+1} + 1$. Now, by Lemma \ref{lemma:rw-corollary}, $\tw(G) \leq 3(R(t, t) -1)2^{2^{s+1}+2} - 1$. This completes the proof. 
\end{proof}

Note that in the above proof, it is sufficient to modify the definition of $\mathcal{F}$ prime decomposition to decompose by 1-cutsets instead of all clique cutsets. 

\end{document}